\documentclass[11pt]{amsart}
\usepackage[english]{babel}
\usepackage{graphics}
\usepackage{epsfig}
\usepackage{latexsym,amssymb,amsfonts,amsmath, amscd, euscript}

\newcommand{\R}{\mathbb{R}}

\newcommand{\Power}{\mathcal P}
\newcommand{\B}{\mathcal B}
\newcommand{\Si}{\mathbb S}

\newtheorem{definition}{{\bf Definition}}[section]
\newtheorem{theorem}[definition]{{\bf Theorem}}

\newtheorem{corollary}[definition]{{\bf Corollary}}
\newtheorem{proposition}[definition]{\noindent {\bf Proposition}}

\newtheorem{claim}[definition]{\noindent {\bf Claim}}

\newtheorem{problem}[definition]{\noindent {\bf Problem}}

\def\proof{{\parindent0pt {\bf Proof.\ }}}
\def\proofref #1 {{\noindent  {\bf Proof} (#1).}\ }
\def\downarrownogap{\downarrow \!\!}
\def\uparrownogap{\uparrow \!\!}

\newtheorem{lemma}[definition]{\noindent {\bf Lemma}}
\def\endproof{\hfill {\kern 6pt\penalty 500
\raise -0pt\hbox{\vrule \vbox to5pt {\hrule width 5pt
\vfill\hrule}\vrule}}}
\def\centerpicture #1 by #2 (#3){\leavevmode
        \vbox to #2{
        \hrule width #1 height 0pt depth 0pt
        \vfill
        \special{pictfile #3}}}
\begin{document}
\title[Retracts of posets]{Retracts of posets: the chain-gap property and the selection property are independent }

\author{Dwight Duffus} 
\address{Department of mathematics, Emory University, Atlanta, Georgia, USA} 
\email{dwight@mathcs.emory.edu}
\author{Claude Laflamme*} 
\address{University of Calgary, Department of Mathematics and Statistics, Calgary, Alberta, Canada T2N 1N4}
\email{laf@math.ucalgary.ca} 
\thanks{*Supported by NSERC of Canada Grant \# 690404} 
\author{Maurice Pouzet**} \address{PCS, Universit\'e Claude-Bernard Lyon1, Domaine de Gerland, B\^at. Recherche [B], 50 avenue Tony-Garnier, F$69365$ Lyon cedex 07, France}
\email{pouzet@univ-lyon1.fr } 
\thanks{**Research completed while the author visited the Mathematics
and Statistics Department of the University of Calgary in July 2006;
the support provided is gratefully aknowledged.}


\date{\today }

\noindent
\begin{abstract}              
Posets which are retract of products of chains are characterized by
means of two properties: \emph{the chain-gap property} and \emph{ the
selection property} (Rival and Wille, 1981 \cite {R-W}). Examples of
posets with the selection property and not the chain-gap property are
easy to find. To date, the Boolean lattice $\Power (\omega_1)/Fin$
was the sole example of lattice without the selection property \cite
{R-W}. We prove that it does not have the chain-gap property. We
provide an example of a lattice which has the chain-gap property but
not the selection property. This answer questions raised in \cite
{R-W}. \end{abstract} \maketitle

\section{Introduction} 
A poset $P$ is a \emph{retract} of a poset $Q$ if there are two
order-preserving maps $f:P\rightarrow Q$ and $g: Q\rightarrow P$ such
that $g\circ f= 1_P$; these maps being respectively called a
\emph{coretraction} and a \emph{retraction}.  The first author and
I.~Rival \cite{D-R} have defined \emph{an order variety} to be a class
of posets closed under direct products and retracts. I.~Rival and
R.~Wille\cite{R-W} characterized members of the order variety
generated by the class of chains as posets satisfying two properties:
\emph{the chain-gap property} and \emph{the selection property}.  And,
they raised the question of their relationship. They gave examples of
lattices with the selection property for which the chain-gap property
fails.  They showed that $\Power (\omega_1)/Fin$, the quotient of
the power set of $\omega_1$ by the ideal $Fin$ of the finite sets,
does not have the selection property. They asked if it has the
gap-property  and we answer this by the negative.

\begin{theorem}\label{nongapfin}
If $E$ is infinite, $\Power (E)/Fin$ does not have the chain-gap property. 
\end{theorem}

They also asked if there is a lattice with the gap property but
without the selection property, and we answer this question
positively.  Our example is a distributive lattice of size $\aleph_1$
which does not embed the ordinal $\omega_1$. It is built from a
Sierpinskization of a subchain $\Si$ of the real line $\R$ which is
$\aleph_1$-dense, that is, $\vert ]a, b[ \cap \Si \vert \geq \aleph_1$
for every $a < b$ in $\Si$, has no end points and has size $\aleph_1$
(the existence of such chains is well-known and easily proved). Let
$\leq_{_{\omega_1}}$ be an ordering on $\Si$ such that the chain $(\Si,
\leq_{\omega_1})$ has order type $\omega_1$, let $\leq_{_\mathbb{R}}$
be the usual ordering on the reals. The Sierpinskization of $\Si$ is the
poset $(\Si, \leq)$ where $\leq$ is the ordering on $\Si$ defined by $x
\leq y$ iff $x \leq_{_{\omega_1}} y$ and $x \leq_{_\mathbb{R}} y$. Let
$L(\Si, \leq)$ be the lattice generated within the lattice of subsets of
$\Si$ by the principal initial segments of $(\Si,\leq)$. So $L(\Si, \leq)$
consists of all the finite unions of finite intersections of sets of
the form $\downarrownogap x$ for $x\in \Si$, where
$\downarrownogap x: =\{ y \in \Si \;\text{and}\; y \leq x\}$. With
this construction in mind we show:

\begin{theorem}\label{main} 
$L(\Si,\leq)$ has the chain-gap property, but not the selection
property.
\end{theorem}

\section{Preliminaries}

If $C$ is a subset of a poset (or quasiorder) $P$, then let
$C^*:=\{x\in P: y\leq x \; \text{ for all } \; y\in C\}$ denote the
set of upper bounds of $C$ and $C_*:= \{x\in P: x\leq y\; \text{ for
all }\; y\in C\}$ denotes the set of lower bounds of $C$. A pair $(A,
B)$ of subsets of $P$ is a \emph{pregap} of $P$ if $A\subseteq B_{*}$
or, equivalently $B\subseteq A^{*}$. A pregap $(A,B)$ is called
\emph{separable} if $A^{*}\cap B_{*}$ is non-empty, otherwise this is
a \emph{gap} of $P$.  We denote by $B(P)$ the set of separable pregaps of
$P$.  Pregaps are quasiordered as follows: $(A, B) \leq (A', B')$ if $A
\leq A'$ and $B' \leq^{^d} B$, where $A \leq A'$ means that every $a
\in A$ is majorized by some $a' \in A'$, and $B' \leq^{^d} B$ means
that every $b \in B$ majorizes some $b' \in B'$.

The {\it cardinality } of a pair $(A, B)$ of subsets of a poset $P$ is
the pair $(\vert A\vert, \vert B\vert)$. We call the pair
\emph{regular} if $\vert A\vert$ and $\vert B\vert $ are both regular,
or one is regular and the other is zero. Say that $(A', B')$ is a
\emph{subpair} of $(A,B)$ if $A'\subseteq A$ and $B'\subseteq B$; if
both pairs are gaps, call $(A',B')$ a \emph{subgap} of $(A,B)$. A gap
$(A,B)$ of $P$ is said to be \emph{minimal} if all subgaps have the
same cardinality as $(A,B)$.  Call a gap $(A, B)$ \emph{irreducible}
if for all subpairs $(A',B')$, $(A',B')$ is a gap if and only if it
has the same cardinality as $(A,B)$. It is straightforward to show
that every gap has a subgap which is minimal. On an other hand,
irreducible gaps are just minimal gaps all of whose subpairs, of its
cardinality, are gaps.

We now come to the main concepts  of this paper.

\begin{definition}
\begin{enumerate}
\item The poset $P$ has the {\it selection property} (the strong
selection property in the terminology of Nevermann, Wille \cite{N-W})
if there is an order-preserving map $\varphi$ from $B(P)$ to $P$ which
associates to every pair $(A, B)\in B(P)$ an element of $A^\ast \cap
B_\ast$.  
\item An order-preserving map $g$ from $P$ into a poset $Q$
\emph{preserves} a gap $(A, B)$ of $P$ if $(g[A],g[B])$ is a gap of
$Q$. If $g$ preserves all gaps of $P$, it is \emph{gap-preserving}. A
poset $Q$ \emph{preserves} a gap $(A, B)$ of $P$ if there is an
order-preserving map $g:P\rightarrow Q$ which preserves $(A,B)$.  The
poset $P$ is said to have the \emph{chain-gap property} if each gap of
$P$ is preserved by some chain.
\end{enumerate}
\end{definition}

The relationship between the chain-gap property and regular
irreducible gaps is given by the following result by Duffus and Pouzet.

\begin{theorem} \cite{D-P} \label{gapproperty}
An ordered set $P$ has the chain-gap property if and only if every gap
of $P$ contains a regular, irreducible gap.
\end{theorem}

In presence of the selection property, they have proved a bit more.

\begin{proposition}\cite{D-P}\label{selection}
Let $(A, B)$ be a minimal gap of $P$ with $\lambda:= \vert A\vert $
and $\mu:= \vert B\vert$ both infinite. If $P$ has the selection
property then there are two chains $C$ and $D$ of type respectively
$cf(\lambda)$ and $cf(\nu)^\ast$ such that $A \leq C$, $D \leq^d B$ and $(C,
D)$ is a gap. \\
%The use of + for ordinal sum is customary in set theory but for ordered sets 
%it is customary to use \oplus for the linear or ordinal sum operation, so it's
%been changed in the next line.  I think that this is the only use in the paper.
Moreover, if $(A, B)$ is an irreducible gap then ordinal sum $C \oplus D$
is a retract of $P$ which preserves $(A,B)$.
\end{proposition}

We illustrate how the above notions relate to a central problem in
the study of retracts of posets, namely to find conditions that a map
$f:P\rightarrow Q$ must satisfy in order to be a coretraction.
Posets $P$ for which maps satisfying these conditions are coretractions
are called \emph{absolute retracts w.r.t. these conditions}.
For example, each coretraction must be an order-embedding. As it is
well-known, \emph {the absolute retracts w.r.t. order-embeddings} are
the complete lattices. Moreover, these are \emph {injective
w.r.t. order-embedding} (that is, every order-preserving map from a
poset $Q$ to $P$ extends to an order-preserving map to every poset
$Q'$ in which $Q$ order-embeds) and there are enough of them in the
sense that every poset order-embeds into a complete lattice, that is, one of them. 
Every coretraction must be gap-preserving. A somewhat similar situation to
the case of order-embeddings was observed by Duffus and Pouzet \cite
{D-P}, and by Nevermann and Rival \cite{N-R}:
\begin{quote}A poset $P$ is an \emph {absolute retract
w.r.t. gap-preserving maps} if and only if it has the selection
property. Moreover, absolute retracts coincide with \emph {injective
objects w.r.t. gap-preserving maps} and there are enough of them.
\end{quote}
The class of absolute retracts is preserved under retraction and
products (Rival and Wille\cite {R-W}), it contains the chains (Duffus,
Rival and Simonovits \cite{D-R-S}) hence the variety generated by the
class of chains. According to Rival and Wille \cite{R-W}:
\begin{quote}A  poset
$P$ embeds by a gap-preserving map into a product of chains iff $P$
has the chain-gap property. 
\end{quote}
The chain-gap property implies that $P$ is a lattice. Every countable
lattice belongs to the variety generated by the class of chains
\cite{P-R}, hence satisfies the chain-gap property.  However, there are many
lattices for which the chain gap property fails (see \cite{D-R},
\cite{R-W}).

We conclude this section with some notation and remarks necessary for
the proof of Theorems \ref{nongapfin} and \ref{main}.

For $E$ any set, let $\Power (E)$ be the Boolean algebra made of all
subsets of $E$ and let $\Power(E)/Fin$ be the quotient of $\Power(E)$
by the ideal $Fin$ of finite subsets of $E$. Define $p:
\Power(E)\rightarrow \Power(E)/Fin$ to be the canonical projection. For $X,
Y \in \Power(E)$, we set $X\leq_{Fin} Y$ if $X\setminus Y\in Fin$;
this defines a quasiorder on $\Power(E)$, its image by $p$ is the
order on $\Power(E)/Fin$.  Since $\Power(E)/Fin$ is a lattice, there
are no gaps of cardinality $(\lambda,\mu)$ where either $\lambda $ or
$\mu$ is finite. Moreover, by a countable diagonalization argument as
first observed by Hadamard \cite{H}, there are no gaps of cardinality
$(\omega , \omega )$ either.
%Wow.  This is due to Hadamard?!?  How did anyone find this out?  I changed
%the wording a bit to be more explicit.  

%To avoid trivial misstatements, we need to restrict to infinite sets at some point.  

To avoid trivialities, let us assume that $E$ is infinite in what follows.
Gaps of $\Power(E)$ under the above quasiorder correspond under $p$
to gaps in the poset $\Power(E)/Fin$, so for notational simplicity all
our discussion regarding gaps in $\Power (E)$ can be translated in the
latter structure if necessary.

We also recall that the usual Hausdorff topology on $\Power(E)$ is
obtained by identifying each subset of $E$ with its characteristic
function and giving the resulting space $\{0,1\}^E$ the product
topology. A basis of open sets consists of subsets of the form
$O(F,G):=\{X\in \Power(E): F\subseteq X \text{ and } G\cap
X=\emptyset \}$, where $F, G$ are finite subsets of $E$. Endowed with
this topology $\Power(E)$ is compact and Hausdorff, therefore  a 
Baire space (i.e. any countable union of closed
sets with empty interior has empty interior).

Now toward Theorem \ref{main}, we recall that in a poset $P$ the
%I used \downarrownogap or \downarrow \!\! to close up the space between
%\downarrow and the set symbol, just as was already done for \downarrow x
%This turned out to be more extensive that I thought, so I may have missed
%closing up some of the spaces.  Also, the subscripts of the arrows get quite
%complicated later -- I hope that I didn't make things unreadable.
\emph{initial segment generated by a subset} $A$ of $P$ is $\downarrownogap
A:= \{x\in P: x\leq y\; \text{for some} \; y\in P\}$. A subset $A$ is
\emph{cofinal} in $P$ if $\downarrownogap A=P$. The \emph{cofinality} of
$P$, $cf(P)$, is the least cardinality of a cofinal subset.  The
notions of \emph{final segment generated by a subset} $A$ and of
\emph{ coinitiality} of $P$ are defined dually and denoted
respectively $\uparrownogap A$ and $ci(P)$. For a singleton $x\in P$, we
use the notation $\downarrownogap x$ instead of $\downarrownogap \{x\}$. If the
reference to $P$ is needed, particularly in case of several orders on
the same ground set, we use the notation $\downarrow_{P} \!\! A$ instead of
$\downarrownogap A$.

\section{Proof of Theorem \ref{nongapfin}} 

Consider $E=T_2$ the binary tree of finite sequences of $0$ and $1$,
and $T_2(n)$ those sequences of length at most $n$.  We denote by
$()$ the empty sequence and by $s.(i)$ the sequence obtained by adding
$i \in \{0,1\}$ to the sequence $s$.  As mentioned above for
notational simplicity we will consider the quasiorder $\leq_{Fin}$ on
$\Power(E)$ as opposed to the poset $\Power(E)/Fin$ itself.

For $B \subseteq \Power(E)$ set $B^c:= \{E\setminus X: X\in B\}$. We will
be particularly interested in the set $\B$ of maximal branches of $T_2$, a
closed subset of $\Power(E)$ with no isolated points. Notice that for
$(A, B)\in \Power(\B)\times \Power(\B)$, $(A,B^c)$ is a pregap if and
only if $A$ and $B$ are disjoint.

\begin{proposition} \label{fsigmagap}
Let $A$ and $B$ be disjoint subsets of $\B$.  Then $(A, B^c)$ is
separable if and only if $A$ and $B$ are covered under inclusion by
disjoint $F_{\sigma}$ sets.
\end{proposition}

\proof If $X \in A^* \cap B^c_*$ separates $(A,B^c)$, then we can
simply let $A' = \bigcup_n \{Y \in \Power(E): Y \setminus X \subseteq
T_2(n) \}$ and $B' = \bigcup_n \{Y \in \Power(E): Y \setminus
(E\setminus X) \subseteq T_2(n) \}$, two disjoint $F_{\sigma}$ sets
covering $A$ and $B$ as required.

Conversely let $A'$, $B'$ be disjoint $F_{\sigma}$ sets covering $A$
and $B \subseteq \B$.  We may assume without loss of generality that
$A'=\bigcup A'_n$ and $B'=\bigcup B'_n$ are increasing chains of
closed sets in $\B$.

\noindent For any fixed $n$, we claim that there must be an integer
$k_n$ such that any $s \in X \cap Y$ has length at most $k_n$ for any
$X \in A'_n$ and $Y \in B'_n$.  Indeed otherwise for infinitely many
$k$ we could find $s_k \in X_k \cap Y_k$ of length at least $k$ for
some $X_k \in A'_n$ and $Y_k \in B'_n$. But then we could find a
subsequence of $\{s_m: m \in {\mathbb N} \}$ converging to a maximal
branch which by closure would be in $A'_n \cap B'_n$, a contradiction.

\noindent We can also assume that the produced sequence $\{ k_n : n
\in {\mathbb N} \}$ is strictly increasing, and we conclude that 
$X=\bigcup_n \{A_n \cap T_2(k_{n+1}) \setminus T_2(k_n)\} \in A^* \cap B^c_*$,
and therefore $(A, B^c)$ is separable.  \endproof

By considering $A$ to be the single branch, and $B=\B \setminus A$,
one concludes that the above result cannot be strengthened to a
covering by disjoint closed sets.

Although the first part of the proof does generalize to any separable
pregap in $\Power(E)$, it is interesting that the converse is not true
as is shown by an example given by Todorcevic \cite{T}. Indeed $A' =\{
\{s \in E: s.0 \in b\}: b \in \B\}$ and $B' =\{ \{s \in E: s.1 \in
b\}: b \in \B \}$ are two disjoint closed sets in $\Power(E)$ forming
a (Luzin) gap in $\Power(E)$.

Since as mentioned above $\Power(E)$ is a Baire space, we further
have:

\begin{corollary} \label{densegap}
If $A \subseteq \B$ and $B=\B \setminus A$ are both dense  then the pair
$(A, B^c)$ is a gap in $\Power(E)$.
\end{corollary}

We finally arrive at the main reason for considering this structure.

\begin{proposition} \label{noregirrgap}
Let $A$ and $B$ be two disjoint subsets of $\B$. If $(A, B^c)$ is a
gap, then it does not contain a regular irreducible gap.
\end{proposition}

\proof 
For $s\in E$ and $D\subseteq \Power(E)$, we set $D(s):= \{X\in
D: s\in X\}$ and $\check D :=\{s\in E: \vert D(s)\vert = \vert D\vert
\}$.

Now observe that for an infinite $D \subseteq \B$, the least element
of $T_2$, namely the empty sequence $()$, belongs to $\check
D$. Moreover if $s\in \check D$, then either $s.(0)$ or $s.(1)$
belongs to $\check D$, so we conclude that $\check D$ contains a
branch and so certainly is infinite.  Moreover, if $\vert D\vert $ is
regular and uncountable, then $\check D$ must contain more than a
chain and is therefore itself is not a chain.

With this, suppose for contradiction that $(A, B^c)$ contains a
regular irreducible gap of size $(\lambda, \mu)$ in $\Power(E)$. This means that
there is a pair $(A', B')$ such that $A'\subseteq A$, $\vert A'\vert
= \lambda$, $B'\subseteq B$, $\vert B'\vert = \mu$ such that $(A',
B'^c)$ is an irreducible gap.  

As noticed in a previous remark, $\lambda $ and $\mu$ must be infinite
and one of them uncountable. With no loss of generality, we may
suppose that this is $\lambda$.  According to the above observation, $\check A'$
is not a chain and $\check B'$ is infinite, hence there are $s\in
\check A'$, $t\in \check B'$ which are incomparable with respect to the
order on $T_2$. Let $A'':= A'(s)$ and $B'':= B'(t)$. We have
$A''\subseteq A'$, $\vert A''\vert = \vert A'\vert = \lambda$,
$B''\subseteq B'$, $\vert B''\vert = \vert B'\vert = \mu$, and therefore 
$(A'',B''^c)$ must be a gap by the irreducibility assumption.
On the other hand for $Z:= \bigcup A''$, we have $X\leq _{Fin}Z \leq_{Fin}Y$ for
every $X\in A''$ and $Y\in B''^c$, a contradiction. \endproof

With this in hand, the proof of Theorem \ref{nongapfin} breaks into two cases.

{\bf Case 1.} $E$ is denumerable.  We deduce Theorem \ref{nongapfin}
as follows. We identify $E$ by $T_2$, and choose $A\subseteq \B$
and $B:= \B \setminus A$ both dense in $\B$. According to
Corollary \ref{densegap}, $(A, B^c)$ is a gap of $\Power(E)$, and according to
Proposition \ref{noregirrgap}, it does not contain a regular irreducible
gap. According to Theorem \ref{gapproperty}, $\Power(E)/Fin$ does not
have the chain-gap property.

{\bf Case 2.}  $E$ is uncountable. Let $E'$ be a denumerable subset of
$E$. The identity map $1_{E'}$ on $E'$ extends to a map $\varphi$ from
$\Power(E')/Fin$ into $\Power(E)/Fin$. This map is
gap-preserving. Thus, if a gap $(A,B)$ in $\Power(E')/Fin$ is not
preserved by a chain, its image $(\varphi[A], \varphi[B])$ cannot be
preserved by a chain. Since $\Power(E')/Fin$ contains such gaps,
$\Power(E)/Fin$ does, too.  Thus, it does not have the chain-gap property.
\endproof

\section{Proof of Theorem \ref{main}}

The proof naturally breaks into two main parts.

\bigskip

\noindent {\bf Part 1:} $L(\Si,\leq)$ does not have the selection property.

\bigskip

\noindent It suffices to prove the following.

\begin{proposition} \label{part1}
\begin{enumerate}
\item[(1)]  $\omega_1$ does not embed into $L(\Si, \leq)$.
\item[(2)]  $L(\Si, \leq)$ has a minimal gap $(A, \emptyset)$ of size $(\aleph_1, 0)$. 
\end{enumerate}
\end{proposition}

Indeed to see how Theorem \ref{main} follows, let $(a_\alpha)_{\alpha
< \omega_1}$ be an enumeration of the elements of $A$. Set $A_\alpha:
= \{a_\beta: \beta < \alpha\}$. If the selection property holds, then
to every pair $(A_\alpha, \emptyset)$ we can associate an element
$x_\alpha \in A^*_\alpha \cap \emptyset_{*}= A^*_\alpha$ such that
$(A_\alpha, \emptyset) \leq (A_{\alpha^\prime}, \emptyset)$ implies
$x_\alpha \leq x_{\alpha^\prime}$. In particular, for $\alpha \leq
\alpha^\prime$ we must have $x_\alpha \leq x_{\alpha^\prime}$. If
$\omega_1$ does not embed into $L(\Si, \leq)$ then the sequence
$x_\alpha$ must be stationary, and in particular has an upper
bound. If $u$ is such an upper bound, then $u \in A^*_\alpha$ for
every $\alpha$, thus $u \in A^*$. This is impossible since $A$ is
unbounded.

\bigskip

\proofref {of Proposition \ref{part1}} We first prove that (2) holds. 

\begin{lemma} \label{parttwo} 

Fix $r \in \Si$ arbitrary and let $A := \{ \downarrow \! \! x: x \in \Si \text{ and } x
\leq_{_\R} r\}$. Then $(A, \emptyset)$ is a minimal gap in $L(\Si,
\leq)$ of size $(\aleph_1, 0)$.
\end{lemma}

\proofref {of Lemma \ref{parttwo}} The proof will follow after these two claims.
%It seems to clarify things a bit to also include Proof (of Lemma   ), etc, to
%help the readers keep track of what statements are being proved.  

\begin{claim} \label{directed}
$(\Si, \leq)$ is up directed.
\end{claim}

\proofref {of Claim \ref{directed}} Let $x, y \in \Si$.  The set $X := \{z\in \Si : z
\leq_{_{\omega_1}}x \text{ or } z \leq_{\omega_1} y\}$ is countable,
but on the other hand the set $Y:= \{z \in \Si: x, y \leq_\R z \}$ is
uncountable.  Thus, $Y\setminus X$ is non empty, and every $z\in Y
\setminus X$ majorizes $x$ and $y$ in $(\Si, \leq)$, proving our
claim.  \endproof

\begin{claim} \label{upperbound}
A subset $B$ of $L(\Si,\leq)$ has an upper bound if and only if
$\bigcup B$ has an upper bound in $(\Si, \leq)$.
\end{claim}

\proofref {of Claim \ref{upperbound}} If $B$ has an upper bound in $L(\Si,\leq)$, then there is some
member $X$ of $L(\Si,\leq)$ which includes every element of $B$, hence
$\bigcup B \subseteq X$. This set $X$ is a finite union of finite
intersections of principal initial segments of $(\Si, \leq)$; hence $B$
is a subset of a finite union $\downarrownogap x_1\cup \cdots \cup
\downarrownogap x_k$ of principal initial segments of $(\Si,\leq)$. Since
$(\Si,\leq)$ is up-directed by Claim \ref{directed}, these is some $x$
which majorizes $x_1, \dots, x_k$, and such an $x$ majorizes $\bigcup
B$. The converse is trivial and the claim is verified.  \endproof

\bigskip

With these claims, the proof of Lemma \ref {parttwo} goes as
follows. From the fact that $\Si$ is $\aleph_1$-dense of size
$\aleph_1$, $A$ has size $\aleph_1$. Next, let's see that  $(A, \emptyset)$ is a gap.

Now if $A$ was bounded, then Claim \ref{upperbound} would imply that
$\bigcup A \subseteq \downarrownogap z$ for some $z\in \Si$. In particular,
the uncountable initial segment of $\Si$ below $r$ under $\leq_\R$
from $\Si$ would be a subset of the countable initial segment of $\Si$
below $z$ under $\leq_{\omega_1}$, a contradiction. $A$ is therefore
unbounded in $L(\Si,\leq)$ and $(A,\emptyset)$ is a gap.

Finally, to show that $(A, \emptyset)$ is a minimal gap amounts to
show that every countable subset $A^\prime$ of $A$ is bounded.  Indeed
let ${\check{A^\prime}}:= \{x \in \Si: \downarrownogap x \in
A^\prime\}$. Then ${\check{A^\prime}}$ is countable thus is bounded in
$\leq_{\omega_1}$. If $y$ is such a bound, then according to Claim
\ref{directed} there is some $z \in \Si$ such that $y \leq z \mbox{ and
} r \leq z$. Then $\downarrownogap z$ is a bound of $A^\prime$.
\endproof

\bigskip

This concludes the verification of statement (2) of the Proposition, and we
now turn our attention to statement (1).

\begin{claim}\label{pair}
Let $a_1, \dots, a_n\in (\Si, \leq )$ and $A: = \downarrownogap a_1\cap
\ldots \cap \downarrownogap a_n$.  Then there are $i, j\leq
n$ such that $A=\downarrownogap a_i\cap
\downarrownogap a_j$.
\end{claim}

\proofref {of Claim \ref{pair}}
Let $i$ such that $a_i\leq_{\omega_1} a_k$ for all $k$, $1\leq k\leq
n$, and let $j$ such that $a_j\leq_{\mathbb{R}} a_k$ for all $k$,
$1\leq k\leq n$.  Then $A = \downarrownogap a_i \cap
\downarrownogap a_j$.  
\endproof
 
 \bigskip
>From this we immediately get:
 
\begin{claim}\label{principalseg}
Every finite intersection of principal initial segments of $(\Si, \leq)$
is of the form $\downarrownogap x \cap \downarrownogap y$
with $x \leq_\R y$ and $y \leq_{\omega_1} x$.
\end{claim}

% Since  $A_\alpha\subseteq \downarrow\hskip -2pt  X_\alpha$, we have $B_\alpha \subseteq (\leftarrow, x_\alpha]$.  If $\alpha\not=0$ then $A_\alpha\not =\emptyset$, hence  $a_\alpha$, the supremum of $A_\alpha$ in $\R$, exists and is majorized by $x_{\alpha}$. 

Now toward a proof that $\omega_1$ does not embed in $L(\Si, \leq)$,
let $(A_\alpha)_{\alpha < \omega_1}$ be an $\omega_1$-sequence of
elements of $L(\Si, \leq)$. According to Claim \ref{principalseg}, for
each $\alpha<\omega_1$, we may write $A_\alpha :=\bigcup \{ A_{\alpha,
i}: i\in I_\alpha \}$ where $I_\alpha$ is a finite set and $A_{\alpha,
i}: = \downarrownogap x_{\alpha,i} \cap \downarrownogap y_{\alpha, i}$
with $x_{\alpha, i} \leq_{_\mathbb{R}} y_{\alpha, i}$ and $y_{\alpha,
i} \leq_{_{\omega_1}} x_{\alpha, i}$.  Set $X_\alpha:= \{x_{\alpha,
i}: i \in I_\alpha \}$, $Y_\alpha:= \{y_{\alpha, i}: i \in I_\alpha\}$
and $Z_{\alpha}:= X_{\alpha}\cup Y_{\alpha}$.

\begin{claim} \label{disjoint}
If $(A_\alpha)_{\alpha < \omega_1}$ is strictly increasing then the
sets $Z_{\alpha}$'s cannot be pairwise disjoint.
\end{claim}

\proofref {of Claim \ref{disjoint}} Let $\alpha<\omega_1$.  Set $x_\alpha =\max_\R (X_\alpha)$,
$\downarrow_\R \! x_\alpha = \{z\in \Si: z\leq_\R x_\alpha\}$, $\downarrow
_\R \!\! A_{\alpha}:= \{z\in \Si: z\leq_{\R} x \mbox{ for some } x \in
A_{\alpha}\}$.  Since $A_\alpha \subseteq A_\beta$ whenever $\alpha
\leq \beta$ we have $\downarrow_\R \!\! A_\alpha \subseteq \downarrow_\R \!\!
A_\beta$.  Since further $\omega_1$ does not embed into the chain $\Si$,
it does not embed into the chain of initial segments of $(\Si,
\leq_\R)$, hence the $\omega_1$-sequence $(\downarrow_\R \!\!
A_\alpha)_{\alpha<\omega_1 }$ is eventually constant.  Let $\alpha_0$
such that $\downarrow_\R \!\!A_\alpha = \downarrow_\R \!\! A_{\alpha_0}$ for
$\alpha_0\leq \alpha <\omega_1$ and define $A := \downarrow_\R \!\!
A_{\alpha_0}$.  Since $A_\alpha \subseteq \   \downarrownogap X_\alpha$, we
have $\downarrow_\R \!\! A_\alpha  \subseteq   \downarrow_\R \!\! x_{\alpha}$.

\noindent Suppose now that that the $Z_{\alpha}$'s are pairwise
disjoint. Then in particular all the $X_{\alpha}$'s are pairwise
disjoint, and therefore there is at most one $\alpha$ such that $A=
\downarrow_\R\!\! x_{\alpha}$ and we may without loss of generality assume
that $x_\alpha \in \Si\setminus A$ for all $\alpha \geq \alpha_0$.
Since $\omega^{*}_1$ does not embed into $\Si$, there is some $x\in \Si
\setminus A$ for which $X_x:=\{\alpha: \alpha_0<\alpha<\omega_1 \mbox{
and } x<_\R x_\alpha \}$ is uncountable.  Since on the other hand the
$Y_\alpha$ are also assumed to be pairwise disjoints, all
$y_{\alpha,i}$'s are therefore distincts, and since $\{ z\in \Si:
z\leq_{\omega_1} x\}$ is countable, there is some $\alpha\in X_x$ such
that $x \leq_{\omega_1} x_\alpha$ and $x \leq_{\omega_1} y_{\alpha,i}$
for all $i \in I_\alpha$. But now for such an $\alpha$ $x <_\R
x_\alpha \leq_\R y_{\alpha,i}$ where $i$ is such that
$x_\alpha=x_{\alpha,i}$, then this implies that $x \in
A_\alpha$. Since $A_\alpha \subseteq \  \downarrow_\R \!\! A_\alpha=A$, we get
$x\in A$ contradicting our definition of $A$.
\endproof

\begin{claim}\label{subchain} 
If there is a strictly increasing $\omega_1$-sequence of elements of
$L(\Si,\leq)$ then there an $\aleph_1$-dense subchain $\Si'$ of $\R$ and a
strictly increasing $\omega_1$-sequence of elements of $L(\Si',\leq')$
for which all $Z_{\alpha}$'s are pairwise disjoint.
\end{claim}
  
\proofref{of Claim \ref{subchain}}
Start with a strictly increasing $\omega_1$-sequence
$(A_\alpha)_{\alpha < \omega_1}$ of members of $L(\Si, \leq)$.  Since
the $Z_{\alpha}$'s are finite, there is an uncountable subset $U$ of
$\omega_1$ and a finite subset $F$ of $\Si$ such that for all $\alpha,
\beta \in U$, $F$ is an initial segment of $Z_\alpha$ w.r.t. $(\Si,
\leq_{\omega_1})$ and $Z_\alpha \cap Z_\beta = F$. That is
$(Z_{\alpha})_{\alpha\in U}$ forms an uncountable $\Delta$-system.
   
Let $x\in \Si$ such that $F \subseteq \downarrow_{\omega_1} \!\! x$ and write $X=
\downarrow_{\omega_1} \!\! x$. Set
$\Si' := \Si \setminus X$, $A^\prime_{\alpha}:=A_{\alpha} \setminus X$,
$A^\prime_{\alpha,i}:=A_{\alpha,i} \setminus X$, $I'_\alpha:= \{i\in
I_{\alpha}: A^\prime _{\alpha,i}\not =
\emptyset\}$.
 
Since $X$ is countable, $\Si'$ is again an $\aleph_1$ dense chain with
no end-points and the well-ordering induced has order type
$\omega_1$. The intersection order $\leq'$ is the order induced by
$\leq$ on $\Si'$.

The $\omega_1$-sequence $(A^\prime_\alpha)_{\alpha < \omega_1}$ is
increasing and since $X$ is countable, it contains a strictly
increasing subsequence $(A^\prime_\alpha)_{\alpha \in U'}$, for some
uncountable subset $U'$ of $U$. Let $\alpha\in U'\setminus
\min(U')$. Then $A^\prime_\alpha\not =
\emptyset$, hence $A^\prime=
\cup\{A'_{\alpha,i}: i\in I'_\alpha\}$. Since $X$ is an initial
segment of $(\Si, \leq )$ it follows that
$A'_{\alpha,i}=\downarrow_{(\Si',\leq')} \!\!  x_{\alpha,i} \  \cap
\downarrow_{(\Si',\leq')}  \!\! y_{\alpha,i}$. Hence $X'_{\alpha}=
\{x_{\alpha,i}: i\in I'_\alpha\}$, $Y'_{\alpha}=\{y_{\alpha,i}: i\in
I_\alpha\}$, and $Z'_\alpha:= X'_{\alpha}\cup Y'_{\alpha}$. Thus the
$Z'_\alpha$ for $\alpha\in U\setminus \min(U)$ are pairwise
disjoint.\endproof

>From Claim \ref{disjoint} and Claim \ref{subchain}, there is no
strictly increasing $\omega_1$-sequence of elements of
$L(\Si,\leq)$. The proof of Proposition \ref{part1} is complete.
\endproof

\bigskip

\noindent {\bf Part 2:}  $L(\Si, \leq)$ has the chain-gap property.

\bigskip

\noindent Let $(A, B)$ be a gap in $L(\Si, \leq)$. 

\begin{lemma}
There is a partition of $L(\Si, \leq )$ into a prime ideal $I$ and a
prime filter $F$ such that $(A, \emptyset ) $ is a gap of $I$ and
$(\emptyset, B)$ is a gap of $F$
\end{lemma}

\proof 
This just follows from the fact that $L(\Si, \leq)$ is a distributive
lattice (see Pouzet, Rival \cite{P-R}).  \endproof

\begin{lemma}
The gap $(\emptyset, B)$ of $F$ can be separated by a chain.
\end{lemma}

\proof 
According to Pouzet-Rival \cite{P-R}, it suffices to show that the
coinitiality of $F$ is countable.

Let $K: = \{x\in \Si: \downarrownogap x \in F\}$. As a a subset of
$\R$, $K$ has a countable coinitiality, so we can select a countable
subset $D$ coinitial in $K$ w.r.t the order $\leq_\R$. Let $U:= \{x
\in K: x \leq_{_{\omega_1}} y \mbox{ for some } y \in D\}$. Since $D$
is countable, $U$ is countable. Moreover, $U$ is coinitial in $K$.
Indeed, let $x \in K$, then there is $x_1 \in D$ such that $x_1
\leq_\R x$.  Now, either $x_1
\leq_{\omega_1} x$, in which case $x_1\leq x$, or $x
<_{\omega_1} x_1$ but then by definition of $U$, $x \in U$. So in
both cases, $x$ majorizes an element of $U$. Let $\check{U}= \{
\downarrownogap x\cap \downarrownogap y: x \in U, y \in U
\}$. Since $U$ is countable, this set is countable. Moreover it is
coinitial in $F$. Indeed, let $a \in F$, then $a$ is of the form $a =
a_1\cup\ldots \cup a_n$ where $a_i :=\downarrownogap x_i \cap
\downarrownogap y_i$. Since $F$ is a prime filter, some $a_i \in
F$ and since $F$ is a filter, $x_i, y_i$ belong to $K$. Now $U$ being
coinitial in $K$ it follows that there are $x^\prime_i, y^\prime_i \in
U$ such that $ x_i^\prime\leq x_i$ and $ y_i^\prime\leq y_i$ hence
$\downarrownogap x_i^\prime\cap \downarrownogap y^\prime\subseteq
a_i\subseteq a$ proving that $\check{U}$ is coinitial in $F$.
\endproof

\begin{lemma}\label{gapof}
The gap $(A, \emptyset)$ in $I$ can be separated by a chain.
\end{lemma}

\proof Elements of $I$ are of the form $a = a_1\cup a_2 \ldots \cup
a_n$ where $a_i = \downarrownogap x_i \cap \downarrownogap y_i$. Since
$I$ is a prime ideal, for every $a_i$ one of the sets
$\downarrownogap x_i, \downarrownogap y_i$ belongs to
$I$. Consequently the set of finite unions of members of $I$ of the
form $\downarrownogap x$ is cofinal in $I$. For a subset $X$ of
$I$ let $\check{X}:= \{x \in \Si: x \in u \; \text{for some} \; u \in X
\}$, in other words $\check{X}= \bigcup X$.

\begin{claim}\label{gapinI}
Let $A$ be a subset of $I$. Then $(A, \emptyset)$ is a gap in $I$ iff
$\Check{A}$ is not contained into a finitely generated initial segment
of $\check{I}$.
\end{claim}

\proof 
If $(A, \emptyset)$ is not a gap in $I$ then from the above
observation there are $x_1, \ldots x_k\in \Si$ such that for
$a:=\downarrownogap x_1\cup \ldots \cup \downarrownogap x_k$ we
have $a\in I$ and $a$ dominates $A$.  This implies $\Check{A}
\subseteq a$, and proves our claim. The converse is obvious.
\endproof

\begin{claim}
Let $A$ be a subset of $I$, then $(A, \emptyset)$ is a gap in $I$ iff
$(\bar{A}, \emptyset )$, where $\bar{A}:= \{ \downarrownogap x
 \text{ such that } \downarrownogap x\subseteq a,
\text{ for some } a \in A\}$, is a gap.
\end{claim}

\proof 
Observe that $\check{A} = \check{\bar{A}}$ and apply Claim \ref{gapinI}.
\endproof

Let $(A, \emptyset)$ be a gap in $I$. We may suppose that $(A,
\emptyset)$ is minimal. Since every elements of $A$ majorize some
element of $\bar{A}$, in order to separate $(A, \emptyset)$ by a chain
it is enough to separate $\bar{A}, \emptyset)$. Let $\leq_{\ast}$ one
of the two orderings $ \leq_{_{\omega_1}}$, $\leq_{_\mathbb{R}}$
restricted to $\check {I}$ and let $\check I_{*}:= (\check I,
\leq_{\ast})$.  We consider two cases:
\begin{enumerate}
\item $\check{A}$ is an unbounded subset of $\check{I}_{*}$ for some $\leq_{\ast}$.
\item $\check{A}$ is a bounded subset of $\check{I}_{*}$  for the two possible  orderings $\leq_{\ast}$.
\end{enumerate}

\noindent{\bf Case (1).}  The ideal $\check{I}$ is unbounded, hence
$\check{I}_{\ast}$ is unbounded too. Let $(c_\alpha)_{\alpha < \mu} \;
(\mu = \omega_1 \; \text{or} \; \omega)$ be an increasing cofinal
sequence of elements of $\check{I}_{\ast}$.  For $\alpha < \mu$, let
$I_\alpha: = \{u \in I$ such that $u
\subseteq\downarrow_{\check{I}_\ast} c_\alpha\}$. Clearly, the
sequence $(I_\alpha)_{\alpha < \mu}$ is increasing. Next $I =
{\underset{\alpha < \mu}{\bigcup}} I_\alpha$. Indeed, let $u \in I$.
There are $x_1, \dots, x_k\in \check{I}$ such that $u
\subseteq\downarrownogap x_1 \cup \ldots \cup \downarrownogap
x_k$. Hence, there is some $c_\alpha$ such that $ x_1, \ldots, x_k
\leq_\ast c_\alpha$. Since for each $i, 1\leq i\leq k$,
$\downarrownogap x_i \subseteq \downarrow_{\check{I}_\ast} x_i
\subseteq \downarrow_{\check{I}_\ast} c_\alpha$ we have $u
\subseteq \downarrow_{\check{I}_\ast} c_\alpha$ thus $u
\in I_\alpha$.  Finally, for every $\alpha$, $A \setminus I_\alpha$ is
non empty. Indeed, since $\check{A}$ is unbounded in $\check{I}_{*}$,
there is some $x \in \check{A}$ such that $c_{\alpha + 1} \leq_\ast
x$. By definition of $\check{A}$, $x \in a$ for some $a \in A$. But
then $a \not\in I_\alpha$.  Pick $a_{\alpha} \in A\setminus I$ for
each $\alpha< \mu$.  Let $A':= \{a_\alpha: \alpha<\mu\}$. Since $(A,
\emptyset)$ is a minimal gap and $\mu$ is regular, $(A', \emptyset)$
is a regular irreducible gap.

\noindent{\bf Case (2)}. 
Note that in this case $\check I$ is not an ideal of $(\Si, \leq)$
(otherwise $\check{A}$ would be bounded in $\check {I}$ and $(A,
\emptyset)$ would no be a gap in $I$).

Since $\check{A}$ is a bounded subset of $\check{I}$
w.r.t. $\omega_1$, it is countable. Hence, there is a least element
$b$ of $(\Si, \leq_{\omega_1})$ for which $\check{A}
\cap\downarrow_{(\Si, \leq_{\omega_1})}b$ contains a subset $B$ which
is not contained into a finitely generated initial segment of
${\check{I}}$. Let $\widetilde B:=\{\downarrownogap x: x\in B\}$. The
pair $(\widetilde B, \emptyset)$ is a gap and in fact a subgap of $(A,
\emptyset)$ too. Hence it suffices to show that $(\widetilde B,
\emptyset)$ is preserved by a chain.

Let $G:=\{z \in \check{I}: B \subseteq \downarrow_{(\Si,
\leq_{\omega_1})} z\}$. Clearly $b$ is a lower-bound of $G$
w.r.t. $\leq_{\omega_1}$. Let $B_1: = B\cap \downarrow_{\Si}G$
and let $B_2 := B \setminus B_1$. Since $B$ is not contained into a
finitely generated initial segment of ${\check{I}}$ there is some
$i\in \{1, 2\}$ such that $B_i$ has the same property.
 
\noindent{\bf Subcase 1}.  $i=2$. In this case, due to the choice of
$b$, $B_2$ is cofinal into $\downarrow_{(\Si,
\leq_{\omega_1})} b$ thus into $B$. Let $B'$ be a cofinal subset of
$B_2$ w.r.t.  $\leq_{\omega_1}$ having order type $\omega$. We claim
that no countable subset of $B^\prime$ can be contained into a
finitely generated initial segment of ${\check{I}}$. Indeed, if there
is one, then there is one, say $B^{\prime\prime}$, which is contained
into some set of the form $\downarrownogap z$, with $z\in \check
I$.  But, since w.r.t. the order $\leq_{\omega_1}$, $B''$ is cofinal
into $B'$, $B'$ is cofinal into $B_2$ and $B_2$ is cofinal into $B$,
$B''$ is cofinal into $B$, hence from $B^{\prime\prime} \subseteq
\downarrow_{(\Si,\omega_1)}z$ we get $z \in G$. With the
fact that $B^{\prime\prime} \subseteq
\downarrow_{(\Si,\leq_{\mathbb{R}})}z$ this implies $B^{\prime\prime}
\subseteq B_1$, contradiction. Let $\widetilde B': = \{
\downarrownogap x: x \in B'\}$.  The property above says that
$(\widetilde B', \emptyset)$ is a gap and every countable subset
too. This gap is regular and irreducible; since $(\widetilde B,
\emptyset)$ contains this gap, it can be preserved by a chain.

\noindent{\bf Subcase 2}. 
Subcase 1 does not hold. Hence $i=1$.  Again, due to the choice of
$b$, $B_1$ is cofinal into $\downarrow_{(\Si,
\leq_{\omega_1})} b$ thus into $B$. Select a cofinal sequence into
$B_1$ with type $\omega$, say $x_0 <_{\omega_1} x_1 <_{\omega_1}
\ldots <_{\omega_1} x_n <_{\omega_1} \ldots$. Observes that $G$ has no
largest element w.r.t. the order $\leq_{\R}$ (otherwise, if $u$ is the
largest element, then we have both $B_1
\subseteq \downarrow_{(\Si, \leq_{\omega_1})} u$ and $B_1
\subseteq \downarrow_{(\Si,\leq_\mathbb{R})}u$, thus $B_1
\subseteq \downarrownogap u$, contradicting the unboundedness of
$B_1$). Hence, the cofinality of $G$ w.r.t. $\leq_{\R}$ is denumerable
and we may select $u_0 <_{\mathbb{R}} u_1 <_\R \ldots <_\R u_n \ldots$
into $G$ forming a cofinal sequence w.r.t. the order $\leq_\R$.

\begin{claim}
There is a sequence $y_0<y_1< \ldots<y_n< \ldots$ of elements of $B_1$
such that $D:= \{y_n: n<\omega\}$ is cofinal in $(B_1,
\leq_{\omega_1})$ and in $(G, \leq_{\R})$.
\end{claim}

\proof 
First, we define $y_0$. Since $x_0<_{\omega_1}b$, $B_1 \bigcap
\downarrow_{(\Si, \leq _{\omega_1})} x_0$ is contained into
a finitely generated initial segment of ${\check{I}}$. Hence $B_1
\bigcap \uparrow_{(\Si, \leq_{\omega_1})}x_0$ is not
contained into a finitely generated initial segment of
${\check{I}}$. In particular,
\begin{equation}\label{notinclusion}
B_1 \bigcap \uparrow_{(\Si, \leq_{\omega_1})}x_0\not
\subseteq \downarrownogap u_0
\end{equation} 
Since $u_0\in G$ we have $B \subseteq \downarrow_{(\Si,
\leq_{\omega_1})} u_0$. From (\ref{notinclusion}) we get
\begin{equation}\label{ineq2}
B_1 \bigcap \uparrow_{(\Si, \leq_{\omega_1})}x_0\not
\subseteq \downarrow_{\leq_\R} u_0 
\end{equation} 
>From (\ref{ineq2}) there is some $y_0\in B_1$ such that $y_0
\geq_{\omega_1} x_0$ and $y_0 \geq_{_\mathbb{R}} u_0$.

Suppose $y_0 < y_1 \ldots < y_n$ be defined with $x_i \leq_{\omega_1}
y_i, u_i \leq_{\R} y_i$. In order to define $y_{n+1}$ select $x_{n_1}$
and $u_{n_1}$ such that:

$$y_n \leq_{\omega_1} x_{n_1}, x_{n+1} <_{\omega_1} x_{n_1} \;
\text{and}\; y_n \leq_{\mathbb{R}} u_{n_1}, u_{n+1} <_{\R}u_{n_1}$$

As above, since $x_{n_1}<_{\omega_1}b$, $B_1 \bigcap \uparrow_{(\Si,
\leq_{\omega_1})}x_{n_1}$ is not contained into a finitely generated
initial segment of ${\check{I}}$ so $B_1 \bigcap \uparrow_{(\Si,
\leq_{\omega_1})}x_{n_1}\not \subseteq \downarrow_{\leq_\R} u_0$ and
thus there is an element, say $y_{n+1}$ such that $x_{n_1}
\leq_{\omega_1} y_{n+1}$ and $u_{n_1} \leq_{_\mathbb{R}}
y_{n+1}$. Clearly, $y_n < y_{n+1}$, $x_{n+1} \leq_{\omega_1} y_{n+1}$
and $u_{n+1}\leq_{\R} y_{n+1}$. From our construction, $D$ is cofinal
in $(B_1, \leq_{\omega_1)}$ and in $(G, \leq_{\R})$.  \endproof

Since $D$ is cofinal in in $(B_1, \leq_{\omega_1})$ and in $(G,
\leq_{\R})$, $D$ is unbounded in $\check I$. But since $\widetilde{D}$
is a chain, it is unbounded in $I$, hence $(\widetilde D, \emptyset)$
is a regular irreducible gap in $I$.  Since $(\widetilde B,
\emptyset)$ contains this gap, it can be preserved by a chain.

With this, the proof of Lemma \ref{gapof} is complete. \endproof

\begin{problem} 
Let $\kappa$ be such that $\omega < \kappa \leq 2^{\aleph_0}$, $\Si$ be
a $\kappa$-dense subchain of $\R$ of size $\kappa$ and $L(\Si, \leq)$
be the distributive lattice associated with a Sierpinskization of
$\Si$. Does $L(\Si, \leq)$ have the chain-gap property?
\end{problem}

\end{document}